\documentclass[12pt,reqno]{amsart}

\usepackage{amsfonts,amsmath,amsthm}
\usepackage{amssymb,epsfig}
\usepackage{enumerate} 
\usepackage{appendix}
\usepackage{comment}

\usepackage{color} 
\definecolor{vert}{rgb}{0,0.6,0}

\usepackage[text={425pt,650pt},centering]{geometry}
\usepackage{graphicx}
\usepackage{epsfig}
\usepackage{tikz,esvect}
\usetikzlibrary{3d,calc,shapes}

\usepackage{caption}
\usepackage{color} 
\definecolor{vert}{rgb}{0,0.6,0}

\numberwithin{figure}{section}

\theoremstyle{plain}
\newtheorem{thm}{Theorem}[section]

\newtheorem{defn}{Definition}

\newtheorem{conjecture}{Conjecture}

\newtheorem{lem}[thm]{Lemma}
\newtheorem{cor}[thm]{Corollary}
\newtheorem{prop}[thm]{Proposition}
\theoremstyle{remark}
\newtheorem{rem}{\bf{Remark}}
\numberwithin{equation}{section}

\newcommand{\N}{\mathbb{N}}

\newcommand{\R}{\mathbb{R}}

\newcommand{\Lip}{{\rm Lip\,}}

\newcommand{\ep}{\varepsilon}

\newcommand{\ol}{\overline}

\begin{document}

\title[Discrete critical Coagulation-Fragmentation equation]
{Discrete Coagulation-Fragmentation equations with multiplicative coagulation kernel and constant fragmentation kernel}

\author[J. JANG, H. V. TRAN]
{Jiwoong Jang, Hung V. Tran}

\thanks{
The work of HT is partially supported by NSF grant DMS-2348305.
}

\address[J. Jang]
{
Department of Mathematics, 
University of Maryland-College Park, William E. Kirwan Hall, 4176 Campus Drive, College Park, Maryland, USA, 20740.}
\email{jjang124@umd.edu}

\address[H. V. Tran]
{
Department of Mathematics, 
University of Wisconsin-Madison, Van Vleck Hall, 480 Lincoln Drive, Madison, Wisconsin, USA, 53706.}
\email{hung@math.wisc.edu}

\keywords{discrete Coagulation-Fragmentation equations; singular Hamilton-Jacobi equations; regularity; well-posedness; Bernstein functions; long-time behaviors; gelation; mass-conserving solutions; viscosity solutions}
\subjclass[2010]{
35B65, 35F21, 35Q99, 44A10, 45J05
}

\maketitle

\begin{abstract}
Here, we study a discrete Coagulation-Fragmentation equation with a multiplicative coagulation kernel and a constant fragmentation kernel, which is critical.
We apply the discrete Bernstein transform to the original Coagulation-Fragmentation equation to get two new singular Hamilton-Jacobi equations and use viscosity solution methods to analyze them.
We obtain well-posedness, regularity, and long-time behaviors of the viscosity solutions to the Hamilton-Jacobi equations in certain ranges, which imply the well-posedness and long-time behaviors of {mass-conserving} solutions to the Coagulation-Fragmentation equation.
The results obtained provide some definitive answers to a conjecture posed in \cite{EMP, ELMP}, and are counterparts to those for the continuous case studied in \cite{Tran-Van-1}.
\end{abstract}


\section{Introduction}

Coagulation-Fragmentation (C-F) equations are integrodifferential equations that arise in various scientific fields, including physics \cite{Smo,BT,VZ}, astronomy \cite{Safro}, and biology, particularly in the study of animal group sizes \cite{Niwa, DLP, LNP}.
Here, coagulation represents the binary merging when two clusters of particles meet and fragmentation represents the binary splitting of a cluster, both at some pre-determined rates. 
With only coagulation and fragmentation governing the dynamics, the C-F describes the evolution of cluster sizes over time.

In this paper, we study a discrete C-F in which cluster sizes are natural numbers.
In this framework, the solution, while still physical, might not conserve mass at all times because of the formation of clusters of infinite size in finite time. This is called a \emph{gelation} phenomenon, which happens when the coagulation is strong enough.
In general, the relative strengths between the coagulation kernel and fragmentation kernel determine whether or not gelation happens regardless of the initial data. 
However, in critical settings where coagulation and fragmentation are of the same relative strength, it is not very clear how solutions behave, and more careful analysis needs to be done based on the initial data, which can be seen below.

\begin{figure}
  \centering
  \begin{minipage}[b]{2.6in}
    \includegraphics[width=2.5in]{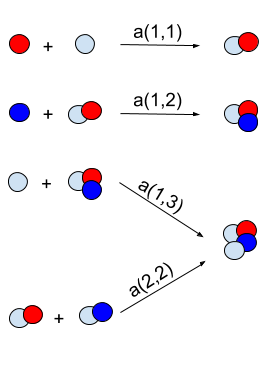}
  \end{minipage}
  \hfill
  \begin{minipage}[b]{2.6in}
    \includegraphics[width=2.5in]{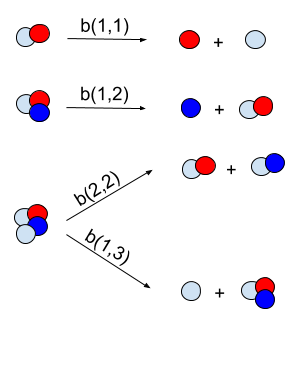}
  \end{minipage}
  \caption{Coagulations and fragmentations}
  \label{fig1}
\end{figure}

Let $\rho(j,t) \geq 0$ be the density of clusters of particles of size $j\in \N$ at time $t\geq 0$.
 We write the discrete Coagulation-Fragmentation equation as
\begin{equation}\label{eq:CF}
\begin{cases}
\rho_t(j,t) = Q_c(\rho)(j,t) + Q_f(\rho)(j,t) \qquad &\text{ on } \N \times (0,\infty),\\
\rho(j,0) = \rho_0(j) \qquad &\text{ on } \N.
\end{cases}
\end{equation} 
Here, the coagulation term $Q_c$ and the fragmentation term $Q_f$ are given by
\[
Q_c(\rho)(j,t)=\frac{1}{2}\sum_{k=1}^{j-1} a(j-k,k) \rho(j-k,t) \rho(k,t) - \rho(j,t) \sum_{k=1}^\infty a(j,k) \rho(k,t),
\]
and
\[
Q_f(\rho)(j,t)=-\frac{1}{2}\rho(j,t)\sum_{k=1}^{j-1} b(j-k,k)  + \sum_{k=1}^\infty b(j,k) \rho(j+k,t).
\]
The coagulation kernel $a(\cdot,\cdot)$ and the fragmentation kernel $b(\cdot,\cdot)$ are nonnegative and symmetric functions on $\N \times \N$.
See Figure \ref{fig1}.
In fact, $Q_c(\rho)$ and $Q_f(\rho)$ can be computed by a simple counting argument.
In this work, we restrict our attention to the critical case of a multiplicative coagulation kernel and a constant fragmentation kernel, that is,
\begin{equation}\label{eq:ab}
a(j,k) = jk \qquad \text{ and } \qquad b(j,k)=1 \qquad \text{ for all } j,k\in \N.
\end{equation}
The initial data $\rho_0:\N \to [0,\infty)$ is given.
This critical discrete equation is the counterpart of the continuous case studied in \cite{EMP, ELMP, Piskunov, Laurencot-2, Tran-Van-1, Mitake-Tran-Van, Tran-Van-2}.

We note that \eqref{eq:CF} is simply a countable system of coupled ordinary differential equations (ODE).
From the first look, the fact that $a(\cdot,\cdot)$ is unbounded makes the analysis of this countable system of ODE intriguing.
For $l\geq 0$, the $l$-th moment of $\rho$ is 
\[
m_l(t)=m_l(\rho)(t)=\sum_{j=1}^\infty j^l \rho(j,t).
\]
When there is no confusion, we use $m_l(t)$ instead of $m_l(\rho)(t)$ for brevity.
In particular, $m_1(t)$ is the total mass of the system at time $t\geq 0$.
Throughout this paper, we also assume
\begin{equation}\label{eq:rho-0}
    m_2(\rho_0) = \sum_{j=1}^\infty j^2 \rho_0(j)<+\infty.
\end{equation}

We refer the reader to \cite{Spouse, Mc, Ball-Carr, Carr, Costa1, Jeon, Laurencot-0, FM, Canizo} for the well-posedness theory and large time behavior of mass-conserving solutions for various discrete C-F.
The two survey papers \cite{Wa, Costa2} give clear and detailed overviews of discrete C-F.
We also refer to the monograph \cite{BLL} and the references therein for recent progress.
It is worth noting that there is a close connection between \eqref{eq:CF} and $3$-wave kinetic equations; see \cite{SoTr} and the references therein.

\subsection{Definition of solutions}
Notice that \eqref{eq:CF} can be formulated as, for any bounded test function $g:\N \to \R$,
\begin{multline}\label{eq:weak form}
    \frac{d}{dt}\sum_{j=1}^\infty g(j) \rho(j,t) = \frac{1}{2}\sum_{j=1}^\infty \sum_{k=1}^\infty \left(g(j+k)-g(j)-g(k) \right) a(j,k) \rho(j,t)\rho(k,t)\\
    -\frac{1}{2}\sum_{j=1}^\infty \sum_{k=1}^{j-1} \left(g(j)-g(k)-g(j-k) \right)b(k,j-k) \rho(j,t).
\end{multline}
Indeed, one can multiply \eqref{eq:CF} by $g(j)$ and sum over $j$ to obtain the above.
In this discrete setting, \eqref{eq:weak form} is equivalent to \eqref{eq:CF}.
For $l\in \N$, let $g(j)=\mathbf{1}_{l}(j)$ be a test function in \eqref{eq:weak form}, where $\mathbf{1}_l$ is the characteristic function of $\{l\}$, to get
\begin{align*}
\rho_t(l,t)=\ &\frac{1}{2}\sum_{j+k=l} a(j,k)\rho(j,t) \rho(k,t) - \rho(l,t)\sum_{k=1}^\infty a(l,k) \rho(k,t)\\
&\ -\frac{1}{2}\sum_{k=1}^{l-1} b(k,l-k) \rho(l,t)+\sum_{j=l+1}^\infty b(l,j-l) \rho(j,t)\\
=\ & Q_c(\rho)(l,t) + Q_f(\rho)(l,t).
\end{align*}

In light of \eqref{eq:ab}, \eqref{eq:weak form} is simplified as
\begin{multline}\label{eq:w}
    \frac{d}{dt}\sum_{j=1}^\infty g(j) \rho(j,t) = \frac{1}{2}\sum_{j=1}^\infty \sum_{k=1}^\infty \left(g(j+k)-g(j)-g(k) \right) jk \rho(j,t)\rho(k,t)\\
    -\frac{1}{2}\sum_{j=1}^\infty \sum_{k=1}^{j-1} \left(g(j)-g(k)-g(j-k) \right)\rho(j,t).
\end{multline}

\begin{defn}
    We say that $\rho$ is a solution to \eqref{eq:CF} with kernels \eqref{eq:ab} if $\rho(\cdot,0)=\rho_0$ and \eqref{eq:w} holds for any bounded test function $g:\N \to \R$.
\end{defn}

\subsection{Some formal calculations and a conjecture}\label{sec:formal cal}
Let $g \equiv 1$ be a test function in \eqref{eq:weak form} to get
\[
m_0'(t) = \frac{1}{2}(m_1(t)-m_1(t)^2) - \frac{1}{2}m_0(t).
\]
Assume that $m_1(t)=m_1(0)=m>0$ as this is true before gelation occurs (if ever).
Then, we get
\begin{equation}\label{eq:m0}
m_0(t)=m-m^2 + e^{-t/2}\left(m_0(0) - (m-m^2)\right).
\end{equation}
Thus, if $m>1$, we see right away that $m_0(t)$ will be negative in finite time, which indicates strongly that gelation occurs.

If $m=1$, then $m_0(t) = e^{-t/2} m_0(0)$, which means that $m_0(t)$ decays exponentially fast to $0$.
Thus, if we do not see gelation, then it will still happen at infinity, and all mass will disappear at infinity.

If $0<m<1$, then $m_0(t)$ converges exponentially fast to $m-m^2$ if there is no gelation.
We thus see a distinctive difference between the case $m=1$ and $0<m<1$.

It is important to note a crucial difference between this discrete C-F and its continuous counterpart.
In the continuous setting, the formation of clusters of zero size, a {\it dust} formation, might happen either at finite or infinite time.
However, in our current discrete setting, there is no dust formation as cluster sizes are always natural numbers.
As a matter of fact, we see in the above formal computations that only gelation occurs in finite time when $m>1$ or possibly at infinity when $m=1$.

The criticality of the C-F can also be seen through this calculation.
Basically, the coagulation term gives $-\frac{1}{2}m_1(t)^2$ and the fragmentation term contributes the leading term $\frac{1}{2}m_1(t)$.
The relative strength of these terms is determined by $m_1(0)=m$ given by the initial data.

\smallskip

The following conjecture was formulated in \cite{EMP, ELMP}.

\begin{conjecture}\label{conj1}
    Consider \eqref{eq:CF} with kernels \eqref{eq:ab} and initial data satisfying \eqref{eq:rho-0}.
    The critical mass of the initial data should be $m=m_1(\rho_0)=1$.
    For $m>1$, there is gelation.
    For $0<m\leq 1$, \eqref{eq:CF} has a unique mass-conserving solution for all time.
\end{conjecture}

See  \cite[Page 163 and Open Problem 2.10]{EMP} as well as \cite[Appendix A]{ELMP}. 
We also refer the reader to \cite[Case (9) on page 261]{VZ}, \cite[Page 65]{Costa2}, and \cite[Section 10.3.4]{BLL}.
Our goal in this paper is to provide some definitive answers to this conjecture.

\subsection{The discrete Bernstein transform and a critical Hamilton-Jacobi equation}
The Bernstein transform was used in the context of C-F in \cite{MP-1, MP-2, MP-3}.
Denote by $F(x,t)$ the discrete Bernstein transform of $\rho(j,t)$, that is, for $(x,t) \in [0,\infty)\times [0,\infty)$,
\[
F(x,t) = \sum_{j=1}^\infty (1-e^{-jx}) \rho(j,t).
\]
For each fixed $x>0$, let $g(j)=1-e^{-jx}$ be a test function in \eqref{eq:w} to get
\begin{align*}
    F_t(x,t)&=- \frac{1}{2}\sum_{j=1}^\infty \sum_{k=1}^\infty \left(1-e^{-jx} \right)\left(1-e^{-kx} \right) jk \rho(j,t)\rho(k,t)\\
    &\quad - \frac{1}{2}\sum_{j=1}^\infty \sum_{k=1}^{j-1} \left(-1-e^{-jx} +e^{-(j-k)x}+e^{-kx} \right)  \rho(j,t)\\
    &=- \frac{1}{2} \left(m_1(t)-F_x(x,t)\right)^2\\
   &\quad -\frac{1}{2}\sum_{j=1}^\infty\left(-j-je^{-jx}+1+e^{-jx}+2e^{-x}\frac{1-e^{(j-1)x}}{1-e^{-x}}\right) \rho(j,t)\\
   &=- \frac{1}{2} \left(m_1(t)-F_x(x,t)\right)^2+\frac{1}{2} \left(m_1(t)+F_x(x,t)\right)-\frac{1}{2}F - \frac{1}{e^x-1}F.
\end{align*}
If we assume that the solution conserves mass for all time, that is, $m_1(t)=m>0$ for all $t\geq 0$ for some given $m>0$, then $F$ solves the following singular Hamilton-Jacobi equation
\begin{equation}\label{eq:HJ}
    \begin{cases}
        F_t +\frac{1}{2}\left(F_x-m\right)\left(F_x-m-1\right)+\frac{1}{2}F + \frac{1}{e^x-1}F -m=0 \qquad &\text { in } (0,\infty)^2,\\
        0\leq F(x,t) \leq m(1-e^{-x}) \qquad &\text { on } [0,\infty)^2,\\
        F(x,0)=F_0(x) \qquad &\text { on } [0,\infty).
    \end{cases}
\end{equation}
In the above, $F_0$ is the discrete Bernstein transform on $\rho_0$. Here, the upper bound of $F$ in the above equation comes from the heuristics
\[
F(x,t) = \sum_{j=1}^\infty (1-e^{-jx}) \rho(j,t) \leq \sum_{j=1}^\infty j(1-e^{-x})\rho(j,t) = m(1-e^{-x}),
\]
with the elementary inequality $1-e^{-jx}\leq j(1-e^{-x})$, which holds only for $j\geq1$.
In the continuous setting, \cite{Tran-Van-1, Mitake-Tran-Van, Tran-Van-2} used the continuous Bernstein transform and the analysis of the corresponding Hamilton-Jacobi equation, which is similar to \eqref{eq:HJ}, to read off information and to give several definitive answers to Conjecture \ref{conj1}.
We develop further this program to study \eqref{eq:CF} with kernels \eqref{eq:ab}. 

Interestingly, in this discrete setting, \eqref{eq:HJ} can be written in another form as follows.
Denote by
\[
G(z,t)=F(-\log z,t) \qquad \text{ for } (z,t)\in (0,1]\times [0,\infty).
\]
Then, for $(z,t)\in (0,1]\times [0,\infty)$,
\[
F_t(-\log z,t)=G_t(z,t), \qquad F_x(-\log z,t) = -zG_z(z,t).
\]
Therefore, \eqref{eq:HJ} can be formulated in an equivalent way as
\begin{equation}\label{eq:G}
    \begin{cases}
        G_t +\frac{1}{2}\left(zG_z+m\right)\left(zG_z+m+1\right)+ \frac{1+z}{2(1-z)}G -m=0 \qquad &\text { in } (0,1)\times(0,\infty),\\
        0\leq G(z,t) \leq m(1-z) \qquad &\text { in } (0,1)\times(0,\infty),\\
        G(z,0)=G_0(z) \qquad &\text { on } (0,1].
    \end{cases}
\end{equation}
Here, $G_0(z)=F_0(-\log z)$ for $z\in (0,1]$.

\smallskip

We need to use both \eqref{eq:HJ} and \eqref{eq:G} in this paper to derive necessary estimates and results.
In the continuous setting \cite{Tran-Van-1}, we only have one singular Hamilton-Jacobi equation similar to \eqref{eq:HJ} where the term $\frac{F}{x}$ appears instead of $\frac{1}{2}F + \frac{1}{e^x-1}F$.
A stark difference here is that the new equation \eqref{eq:G}, which does not appear in the continuous case, plays an essential role in our analysis, and we need to rely crucially on both \eqref{eq:HJ} and \eqref{eq:G}.
To the best of our knowledge, \eqref{eq:HJ} and \eqref{eq:G} have not appeared in the literature.

\subsection{Main results}
We first give our main results of the existence/nonexistence of mass-conserving solutions to \eqref{eq:CF}.

\smallskip

We consider the case $m_1(0)>1$ and show that there is no mass-conserving solution to \eqref{eq:CF} for all time.

\begin{thm}\label{thm:nonexistence-of-C-F}
    Assume the initial data $\rho_0:\mathbb{N}\to[0,\infty)$ satisfies $m_1(0)>1$. 
    Assume the coagulation and fragmentation kernels are given as in \eqref{eq:ab}.
    Then, there is no mass-conserving solution to \eqref{eq:CF} for all time.
    \end{thm}

Next, we focus on the case $0<m_1(0) \leq 1$ and obtain the uniqueness result for mass-conserving solutions to \eqref{eq:CF} in the following proposition.

\begin{prop}\label{prop:uniqueness-of-C-F}
    Assume the initial data $\rho_0:\mathbb{N}\to[0,\infty)$ satisfies $m_1(0)=m\in(0,1]$. 
    Assume the coagulation and fragmentation kernels are given as in \eqref{eq:ab}. 
    Then, there exists at most one mass-conserving solution $\rho:\mathbb{N}\times[0,\infty)\to[0,\infty)$ to \eqref{eq:CF} with $\rho(\cdot,0)=\rho_0$.
\end{prop}

Theorem \ref{thm:nonexistence-of-C-F} and Proposition \ref{prop:uniqueness-of-C-F} are counterparts to those for the continuous case studied in \cite{Tran-Van-1}.
Their proofs are similar to those of \cite[Corollaries 1.3, 1.5]{Tran-Van-1}.

It is a much harder task to show the existence of mass-conserving solutions to \eqref{eq:CF}.
We prove this existence result under the assumptions that $m_1(0)\in\left(0,\frac{1}{2}\right)$ and \eqref{eq:rho-0} holds.

\begin{thm}\label{thm:existence-of-C-F}
    Assume the initial data $\rho_0:\mathbb{N}\to[0,\infty)$ satisfies $m_1(0)=m\in\left(0,\frac{1}{2}\right)$ and \eqref{eq:rho-0}. Assume the coagulation and fragmentation kernels are given as in \eqref{eq:ab}. Then, there exists a unique mass-conserving solution $\rho:\mathbb{N}\times[0,\infty)\to[0,\infty)$ to \eqref{eq:CF} with $\rho(\cdot,0)=\rho_0$.
    \end{thm}

Basically, to obtain Theorem \ref{thm:existence-of-C-F}, we need to show that $F$ and $G$, which are solutions to \eqref{eq:HJ} and \eqref{eq:G} respectively, are smooth and $\partial^k_z G \leq 0$ in $(0,\infty)^2$ for all $k\in \N$.
To do so, we need to study the regularity of the viscosity solutions to \eqref{eq:HJ} and \eqref{eq:G} and utilize their specific structures.
In the continuous setting \cite{Tran-Van-1}, it was enough to show that $F$ is smooth and is a Bernstein function, that is, $(-1)^{k+1}\partial^k_x F \geq 0$ in $(0,1)\times(0,\infty)$ for all $k\in \N$.
In this discrete setting, we cannot deal with $F$ directly, and we instead show that $\partial^k_z G \leq 0$ in $(0,1)\times(0,\infty)$ for all $k\in \N$, which is a stronger requirement.
It is important to note that the unique mass-conserving solution obtained in Theorem \ref{thm:existence-of-C-F} is a solution  to \eqref{eq:CF} in the classical sense.
Thus, Theorem \ref{thm:existence-of-C-F} is more desirable than \cite[Theorem 1.8]{Tran-Van-1} where only weak solutions in the measure sense were considered.
This highlights another crucial difference between the discrete and continuous situations.

\smallskip

We next present the large-time behavior result of mass-conserving solutions to \eqref{eq:CF} under the setting of Theorem \ref{thm:existence-of-C-F}.

\begin{thm}\label{thm:m3}
    Under the same setting of Theorem \ref{thm:existence-of-C-F} with the same notations, the density $\rho(l,t)$ converges to $\tilde{\rho}(l)$ as $t\to\infty$ for each $l\in\mathbb{N}$. Moreover, the equilibrium $\tilde{\rho}:\mathbb{N}\to[0,\infty)$ is a unique stationary solution to \eqref{eq:CF} with $m_1(\tilde{\rho})=m$, characterized by the recursive formula, for $l\in \N$,
    \begin{equation}\label{eq:recursive}
    \tilde{\rho}(l)=\frac{2m(1-m)}{(2m+1)l+1}+ \frac{1}{(2m+1)l+1} \left(\sum_{i=1}^{l-1}i(l-i) \tilde{\rho}(i)\tilde{\rho}(l-i)-2 \sum_{i=1}^{l-1}\tilde{\rho}(i) \right).
    \end{equation}
    \end{thm}

For large-time behavior results and self-similar solutions to C-F, we refer the reader to \cite{MP-1, MP-2, MP-3, NV, DLP, LNP,  Laurencot-1, CFFV} and the references therein for the continuous setting and \cite{Carr, FM, Canizo} and the references therein for the discrete setting.

\smallskip

In the above theorem, under the assumptions $m_1(0)=m\in\left(0,\frac{1}{2}\right)$ and \eqref{eq:rho-0}, we prove that the mass-conserving solution $\rho(\cdot,t)$ converges to the unique non-trivial stationary solution $\tilde \rho$ with $m_1(\tilde{\rho})=m$ as $t\to \infty$, which satisfies the recursive formula \eqref{eq:recursive}.
This phenomenon is completely different from that of the continuous case in \cite[Theorem 1.9]{Tran-Van-1}, where $\rho(\cdot,t)$ turns to dust as $t\to \infty$.

\smallskip

On the other hand, for $m \geq 1$, \eqref{eq:CF} with kernels given as in \eqref{eq:ab} does not admit a stationary solution $\tilde \rho:\N \to [0,\infty)$ such that $m_1(\tilde{\rho})=m$ (see Lemma \ref{lem:series-m-big}).
Of course, the case $m>1$ is consistent with the result of Theorem \ref{thm:nonexistence-of-C-F}.
The case $m=1$ in Lemma \ref{lem:series-m-big} implies that gelation happens and all the mass is gone at infinity, which agrees with the heuristic discussion in Section \ref{sec:formal cal}.
Here, the situation differs significantly from the continuous case, where a one-parameter family of stationary solutions with unit mass exists \cite{Mitake-Tran-Van}.

For $m\in (0,1/2]$, Lemma \ref{lem:series-m-0-1/2} confirms that \eqref{eq:CF} with kernels given as in \eqref{eq:ab} admits a unique stationary solution $\tilde \rho:\N \to [0,\infty)$ such that $m_1(\tilde{\rho})=m$, which satisfies \eqref{eq:recursive}.
Here, $m=1/2$ is included because of the continuity with respect to $m$ of $\tilde \rho(l)=\tilde \rho_m(l)$ in \eqref{eq:recursive}.
Through various numerical computations, we believe that \eqref{eq:recursive} should give us the unique stationary solution for the range $1/2<m<1$ too.
This motivates us to pose the following conjecture.

\begin{conjecture}
    Consider \eqref{eq:CF} with kernels \eqref{eq:ab} and initial data satisfying \eqref{eq:rho-0}.
    \begin{itemize}
        \item[(i)] For $1/2<m<1$, \eqref{eq:CF} admits a unique stationary solution $\tilde \rho:\N \to [0,\infty)$ such that $m_1(\tilde{\rho})=m$. Further, $\tilde \rho$ satisfies the recursive formula \eqref{eq:recursive}.

        \item[(ii)] For $1/2 \leq m_1(\rho_0) <1$, \eqref{eq:CF} admits a unique mass-conserving solution $\rho$ for all time. This solution $\rho$ is in fact a classical solution.
    \end{itemize}
\end{conjecture}

\subsection*{Organization of the paper}
In Section \ref{sec:Well-posedness-of-F}, we prove Theorem \ref{thm:nonexistence-of-C-F} and Proposition \ref{prop:uniqueness-of-C-F} by investigating the well-posedness of viscosity solutions to \eqref{eq:HJ}. In Section \ref{sec:regularity-F}, we establish regularity results of solutions $F$ to \eqref{eq:HJ}, and in Section \ref{sec:solution-CF}, we prove Theorem \ref{thm:existence-of-C-F} by constructing solutions to \eqref{eq:CF} from regularity results of solutions $G$ to \eqref{eq:G} under the assumptions $m_1(\rho_0)\in\left(0,\frac{1}{2}\right)$ and \eqref{eq:rho-0} on the initial data $\rho_0$. We prove Theorem \ref{thm:m3} in Section \ref{sec:large-time-behavior} under the same assumptions.

\section{Well-posedness of viscosity solutions to \eqref{eq:HJ}}\label{sec:Well-posedness-of-F}
In this section, we study the well-posedness of viscosity solutions to \eqref{eq:HJ}.

\subsection{The case of $m\in (0,1]$}
We show in this case that \eqref{eq:HJ} is well-posed in the following sense.

\begin{prop}\label{prop:m1}
    Assume $m\in (0,1]$.
    Let $F_0 \in \Lip([0,\infty))$ be such that $0 \leq F_0(x) \leq m(1-e^{-x})$ for $x \in [0,\infty)$.
    Then \eqref{eq:HJ} has a unique Lipschitz viscosity solution $F$.
    \end{prop}

As the coefficient $1/(e^x-1)$ of the zeroth order term in \eqref{eq:HJ} is singular at $0$, we cut off the singularity by introducing, for $n\in \N$,
\[
  \theta_n(x) = \max\left\{\frac{1}{n},e^x-1 \right\} \qquad \text{ for } x \in [0,\infty).
\]
For each $n\in \N$, we consider the following approximation of \eqref{eq:HJ}
\begin{equation}\label{eq:HJ-n}
    \begin{cases}
        F_t +\frac{1}{2}\left(F_x-m\right)\left(F_x-m-1\right)+\frac{1}{2}F + \frac{1}{  \theta_n(x)}F -m=0 \qquad &\text { in } (0,\infty)^2,\\
        F(0,t)=0 \qquad &\text { on } [0,\infty),\\
        F(x,0)=F_0(x) \qquad &\text { on } [0,\infty).
    \end{cases}
\end{equation}
\begin{lem}\label{lem:F-n}
    Fix $m\in (0,1]$.
    For each $n\in \N$, \eqref{eq:HJ-n} has a unique bounded viscosity solution $F^n$.
    Moreover, $\{F^n\}$ is equi-Lipschitz on $[0,\infty)^2$, and $F^{n+1} \leq F^n$ for all $n\in \N$.
\end{lem}
\begin{proof}
By the classical Perron method \cite{Tran}, \eqref{eq:HJ-n} has a unique solution $F^n$ satisfying
\begin{equation}\label{eq:Fn-1}
F_0(x) - Ct \leq F^n(x,t) \leq F_0(x)+Ct \qquad \text{ for } (x,t) \in [0,\infty)^2,
\end{equation}
where
\[
C=\frac{1}{2}\|\left((F_0)_x-m\right)\left((F_0)_x-m-1\right)\|_{L^\infty([0,\infty))}+3m.
\]
Here, we used the fact that $0\leq F_0(x) \leq m(1-e^{-x})$, and in particular,
\[
0\leq \frac{F_0(x)}{  \theta_n(x)}\leq \frac{m(1-e^{-x})}{  \theta_n(x)} \leq m.
\]
Besides, we note that $\eta^- \equiv 0$ is a subsolution to \eqref{eq:HJ-n} and $\eta^+ \equiv m$ is a supersolution to \eqref{eq:HJ-n}.
To see the subsolution property, we check that
\[
\frac{1}{2}m(m+1) - m = \frac{1}{2}m(m-1) \leq 0.
\]
Moreover, $\zeta^+(x,t)=m(1-e^{-x})+1/n$ is another supersolution to \eqref{eq:HJ-n} as
\begin{align*}
\frac{1}{2}(me^{-x}-m)(me^{-x}-m-1)+\frac{1}{2}&\left(m(1-e^{-x})+\frac{1}{n}\right) + \frac{m(1-e^{-x})+\frac{1}{n}}{  \theta_n(x)} - m\\
&\geq
\begin{cases}
    \frac{1}{2}m^2(1-e^{-x})^2, \qquad &x \geq \log(1+\frac{1}{n}),\\
    1-m, \qquad &x \leq \log(1+\frac{1}{n}),
\end{cases}
\end{align*}
which is always nonnegative. We thus get,  for $(x,t) \in [0,\infty)^2$,
\begin{equation}\label{eq:Fn-2}
\max\left\{0,F_0(x) - Ct\right\} \leq F^n(x,t) \leq \min\left\{m(1-e^{-x})+\frac{1}{n},F_0(x)+Ct\right\}.
\end{equation}

It is clear that \eqref{eq:Fn-1} (or \eqref{eq:Fn-2}) implies $\|F^n_t\|_{L^\infty([0,\infty)^2)} \leq C$.
We use this in \eqref{eq:HJ-n} to get further that
\begin{equation}\label{eq:Fn-3}
    \|F^n_t\|_{L^\infty([0,\infty)^2)}+ \|F^n_x\|_{L^\infty([0,\infty)^2)} \leq C.
\end{equation}

Next, we use the fact that $  \theta_n \geq \theta_{n+1}$ to yield that $F^{n+1} \leq F^n$.
The proof is complete.
\end{proof}

We are now ready to prove Proposition \ref{prop:m1}.

\begin{proof}[Proof of Proposition \ref{prop:m1}]
    We first show the existence of bounded Lipschitz viscosity solutions to \eqref{eq:HJ}.
    Thanks to Lemma \ref{lem:F-n}, $F^n \to F$ locally uniformly on $[0,\infty)^2$, where $F \in \Lip([0,\infty)^2)$, and for $(x,t) \in [0,\infty)^2$,
\begin{equation*}
\max\left\{0,F_0(x) - Ct\right\} \leq F(x,t) \leq \min\left\{m(1-e^{-x}),F_0(x)+Ct\right\}.
\end{equation*}
By the stability of viscosity solutions \cite{Tran}, we see that $F$ is a  bounded Lipschitz viscosity solution to \eqref{eq:HJ}.
We note that $F$ is also a subsolution to \eqref{eq:HJ-n} for every $n\in \N$.

We now prove the uniqueness of bounded Lipschitz viscosity solutions to \eqref{eq:HJ}.
Let $\tilde F$ be another solution to \eqref{eq:HJ}.
We claim that $\tilde F^n = \tilde F+1/n$ is a supersolution to \eqref{eq:HJ-n} for every $n\in \N$.
To check this, we only need to consider \eqref{eq:HJ-n} in $(0,\log(1+1/n)) \times (0,\infty)$.
For $(x,t) \in (0,\log(1+1/n)) \times (0,\infty)$, $  \theta_n(x)=1/n$, and thus,
\[
\frac{1}{2}\tilde F^n + \frac{1}{  \theta_n(x)}\tilde F^n \geq \frac{1}{2}\tilde F +1 
\geq \frac{1}{2}\tilde F +m
\geq \frac{1}{2}\tilde F +\frac{1}{e^x-1}\tilde F.
\]
Thus, $\tilde F^n $ is a supersolution to \eqref{eq:HJ-n} for every $n\in \N$.
By the usual comparison principle for \eqref{eq:HJ-n} (see \cite{Tran}), $\tilde F^n \geq F$.
Let $n\to\infty$ to imply $\tilde F \geq F$.
By a symmetric argument, we obtain $\tilde F=F$.
The proof is complete.
\end{proof}

It is important to note that we use crucially the property that $m \in (0,1]$ in the claim that $\eta^-\equiv 0$ is a subsolution to \eqref{eq:HJ-n}.

\begin{proof}[Proof of Proposition \ref{prop:uniqueness-of-C-F}]
    Assume that $\rho^1, \rho^2$ are two mass-conserving solutions of \eqref{eq:CF}.
    Let $F^1, F^2$ be the discrete Bernstein transform of $\rho^1, \rho^2$, respectively.

    Thanks to Proposition \ref{prop:m1}, we have $F^1=F^2$.
    Therefore, $\rho^1=\rho^2$, which concludes the proof.
\end{proof}
\subsection{The case of $m\in (1,\infty)$}
We show in this case that \eqref{eq:HJ} does not admit a $C^1$ solution for all time as in the following.

\begin{prop}\label{prop:m2}
    Assume $m>1$.
    Let $F_0 \in C^\infty([0,\infty))$ be such that $0 \leq F_0(x) \leq m(1-e^{-x})$ for $x \in [0,\infty)$.
    Then \eqref{eq:HJ} does not admit a solution $F\in C^1([0,\infty)^2)$.
    \end{prop}

\begin{proof}[Proof of Proposition \ref{prop:m2}]
We proceed by contradiction and assume otherwise that \eqref{eq:HJ} admits a solution $F\in C^1([0,\infty)^2)$.
It is clear that, for $t\geq 0$,
\[
F(0,t)=F_t(0,t)=0.
\]
Let $x\to 0^+$ in \eqref{eq:HJ} and use the point that
\[
\lim_{x\to 0^+}\frac{F(x,t)}{e^x-1} =\lim_{x\to 0^+}\frac{F(x,t)-F(0,t)}{x}\frac{x}{e^x-1} = F_x(0,t)
\]
to imply
\[
\frac{1}{2}(F_x(0,t)-m)(F_x(0,t)-m-1) + F_x(0,t)-m=0.
\]
Hence, $F_x(0,t)=m$ or $F_x(0,t)=m-1$.
Either ways, $F_x(0,t) \geq m-1>0$.
Fix $\sigma \in (0,m-1)$.
For each $t>0$, there exists $x_t \in (0,\infty)$ such that
\[
\varphi(t)= \max_{x\in [0,\infty)} \left( F(x,t) - \sigma x \right) = F(x_t,t) - \sigma {x_t}>0.
\]
It is clear that $F_x(x_t,t)=\sigma$.
Plug this into \eqref{eq:HJ} to get
\begin{align*}
    F_t(x_t,t) &=-\frac{1}{2}\left(F_x(x_t,t)-m\right)\left(F_x(x_t,t)-m-1\right)-\frac{F(x_t,t)}{2} - \frac{F(x_t,t)}{e^{x_t}-1} +m\\
    &\leq -\frac{1}{2}\left(\sigma-m\right)\left(\sigma-m-1\right)-\frac{\sigma x_t}{2} - \frac{\sigma x_t}{e^{x_t}-1} +m\\
    &\leq -\frac{1}{2}\left(\sigma-m\right)\left(\sigma-m-1\right)-\sigma +m = -\frac{1}{2}\left(\sigma-m\right)\left(\sigma-(m-1)\right) =: \delta<0.
\end{align*}
Here, we used $\frac{1}{2}+\frac{x}{e^x-1}\geq1$ for $x>0$ in the second inequality. Therefore, for $t>0$,
\begin{align*}
    \varphi'(t) &= \lim_{s\to 0^+} \frac{\varphi(t)-\varphi(t-s)}{s}\\
    &=\lim_{s\to 0^+}\frac{(F(x_t,t)-\sigma x_t) - (F(x_{t-s},t-s) - \sigma x_{t-s})}{s}\\
    &\leq \lim_{s\to 0^+}\frac{(F(x_t,t)-\sigma x_t) - (F(x_{t},t-s) - \sigma x_{t})}{s}
    =F_t(x_t,t) \leq \delta<0.
\end{align*}
Thus, we have that in the viscosity sense, 
\[
\varphi'(t) \leq \delta<0 \qquad \text{ for } t>0.
\]
Then, for $T=1+\varphi(0)/|\delta|$, we see that $\varphi(T)<0$, which is absurd.
\end{proof}

\medskip

\begin{proof}[Proof of Theorem \ref{thm:nonexistence-of-C-F}]
Assume by contradiction that \eqref{eq:CF} has a mass-conserving solution $\rho:\N \times [0,\infty) \to \R$.
Let $F$ be the discrete Bernstein transform of $\rho$.
Then, $F \in C^1([0,\infty)^2) \cap C^\infty((0,\infty)^2)$ is a viscosity solution to \eqref{eq:HJ}.
Proposition \ref{prop:m2} gives a contradiction immediately.
\end{proof}
\section{Regularity of solutions to \eqref{eq:HJ}}\label{sec:regularity-F}
To study the regularity of the viscosity solution to \eqref{eq:HJ}, we consider the vanishing viscosity process
    \begin{equation}\label{eq:add-ep}
    \begin{cases}
        F_t +\frac{1}{2}\left(F_x-m\right)\left(F_x-m-1\right)+\frac{F}{2} + \frac{F}{e^x-1} -m=\ep a(x) F_{xx} \qquad &\text { in } (0,\infty)^2,\\
        F(0,t)=0 \qquad &\text { on } [0,\infty),\\
        F(x,0)=F_0(x) \qquad &\text { on } [0,\infty),
    \end{cases}
\end{equation}
Here, $a\geq 0$ will be chosen appropriately.

\subsection{Spatial concavity of $F$ when $m\in (0,1]$}
In this subsection, we always choose $a\equiv 1$.
\begin{lem}\label{lem:F-ep-increase}
    Assume $m\in (0,1]$, $a\equiv 1$, and \eqref{eq:rho-0}.
    For $\ep\in (0,1)$, let $F^\ep$ be the unique solution to \eqref{eq:add-ep}.
    Then, 
    \[
    0\leq F^\ep_x \leq m \qquad \text{ on } [0,\infty)^2.
    \]
\end{lem}

\begin{proof}
We first note that
\[
0\leq F^\ep(x,t) \leq m(1-e^{-x})  \qquad \text{ for } (x,t)\in [0,\infty)^2,
\]
which gives
\begin{equation}\label{eq:F-ep-1}
    0\leq F^\ep_x(0,t) \leq m \qquad \text{ for } t \in [0,\infty).
\end{equation}
By the classical regularity for parabolic equations, we have that $F^\ep \in C^\infty((0,\infty)^2) \cap C^2_1([0,\infty)^2)$.
Here, $C^2_1([0,\infty)^2)$ is the space of functions which are $C^2$ in $x$
and $C^1$ in $t$ on $[0,\infty)^2$.

Firstly, for $h>0$ fixed, denote by $\tilde F^\ep(x,t)=F^\ep(x+h,t)$ for $(x,t)\in [0,\infty)^2$.
It is clear that
\[
\frac{\tilde F^\ep(x,t)}{e^{x+h}-1} \leq \frac{\tilde F^\ep(x,t)}{e^x-1},
\]
which implies that $\tilde F^\ep$ is a supersolution to \eqref{eq:add-ep}.
By the usual comparison principle, $\tilde F^\ep \geq F^\ep$.
Therefore,
\[
F^\ep_x \geq 0 \qquad \text{ on } [0,\infty)^2.
\]

Next, we prove that $F^\ep_x \leq m$.
We note that $\tilde F^\ep$ is the bounded, Lipschitz viscosity solution to 
\begin{equation*}
    \begin{cases}
        F_t +\frac{1}{2}\left(F_x-m\right)\left(F_x-m-1\right)+\frac{F}{2} + \frac{F}{e^{x+h}-1} -m=\ep F_{xx} \qquad &\text { in } (0,\infty)^2,\\
        F(0,t)=F^\ep(h,t) \qquad &\text { on } [0,\infty),\\
        F(x,0)=F_0(x+h) \qquad &\text { on } [0,\infty),
    \end{cases}
\end{equation*}
Set $\bar F^\ep(x,t)=F^\ep(x,t)+mh$ for $(x,t)\in [0,\infty)^2$. 
We claim that $\bar F^\ep$ is a supersolution to the above.
To check this, we only need to show that
\[
\frac{F^\ep(x,t)}{e^x-1} \leq \frac{F^\ep(x,t)+mh}{e^{x+h}-1} + \frac{mh}{2},
\]
which follows from the inequality 
\[
mh\left(\frac{1}{2}+\frac{1}{e^{x+h}-1}\right) \geq \frac{mx(e^{x+h}-e^x)}{(e^x-1)(e^{x+h}-1)} \geq \frac{F^\ep(x,t)}{e^x-1} - \frac{F^\ep(x,t)}{e^{x+h}-1}.
\]
Thus, $\tilde F^\ep \leq \bar F^\ep$, which gives
\[
F^\ep_x \leq m \qquad \text{ on } [0,\infty)^2.
\]
The proof is complete.
\end{proof}

\begin{lem}\label{lem:F-ep-concave}
    Assume $m\in (0,1]$, $a\equiv 1$, and \eqref{eq:rho-0}.
    For $\ep\in (0,1)$, let $F^\ep$ be the unique solution to \eqref{eq:add-ep}.
    Then, 
    \[
     F^\ep_{xx} \leq 0 \qquad \text{ on } [0,\infty)^2.
    \]
\end{lem}

\begin{proof}
Write $\phi=F^\ep_{xx}$.
We note that at $F^\ep(0,t)=0$, which gives $F^\ep_t(0,t)=0$ for all $t\geq 0$.
Combine this with the fact that $0\leq F^\ep_x(0,t) \leq m\leq 1$ to yield
\begin{equation}\label{eq:F-ep-xx-0}
    \phi(0,t) = F^\ep_{xx}(0,t) = \frac{1}{2\ep}(F^\ep_x(0,t)-m)(F^\ep_x(0,t)-m+1) \leq 0.
\end{equation}

Differentiate the equation in \eqref{eq:add-ep} with respect to $x$ twice  to get
\begin{multline}\label{eq:F-ep-xx}
    \phi_t +F^\ep_x  \phi_x -\left(m+\frac{1}{2}\right) \phi_x - \ep \phi_{xx} +\phi^2\\
    +\frac{\phi}{2} +\frac{\phi}{e^x-1}  +\frac{e^x}{(e^x-1)^2}\left(\frac{e^x+1}{e^x-1}F^\ep-2F^\ep_x\right)=0.
\end{multline}
By Taylor's expansion, for each $(x,t) \in (0,\infty)^2$, there exists $\theta=\theta(x,t)\in (0,1)$ such that
\[
0=F^\ep(0,t)=F^\ep(x,t) -xF^\ep_x(x,t) +\frac{x^2}{2}F^\ep_{xx}(\theta x,t).
\]
Plug this into \eqref{eq:F-ep-xx} to yield
\begin{multline*}
    \phi_t +F^\ep_x  \phi_x -\left(m+\frac{1}{2}\right) \phi_x - \ep \phi_{xx} +\phi^2\\
    +\frac{\phi}{2} +\frac{\phi}{e^x-1}-\frac{x^2e^x(e^x+1)}{2(e^x-1)^3}F^\ep_{xx}(\theta x,t)  +\frac{e^x}{(e^x-1)^2}\left(x\frac{e^x+1}{e^x-1}-2\right)F^\ep_x=0.
\end{multline*}
We use the fact that $x\frac{e^x+1}{e^x-1}-2\geq 0$ and $0\leq F^\ep_x(x,t)\leq m$ for all $(x,t)\in [0,\infty)^2$ to imply further that
\begin{multline}\label{eq:F-ep-xx-r}
    \phi_t +F^\ep_x  \phi_x -\left(m+\frac{1}{2}\right) \phi_x - \ep \phi_{xx} +\phi^2\\
    +\frac{\phi}{2} +\frac{\phi}{e^x-1}-\frac{x^2e^x(e^x+1)}{2(e^x-1)^3}F^\ep_{xx}(\theta x,t)  \leq 0.
\end{multline}
We now use the classical maximum principle to conclude that $\phi \leq 0$.
Assume otherwise that for some $T>0$, there exists $x_0\geq 0$ such that
\[
\max_{[0,\infty)\times [0,T]} \phi = \phi(x_0,T)>0.
\]
In light of \eqref{eq:F-ep-xx-0}, $x_0>0$.
By the usual maximum principle and the point that $F^\ep_{xx}(\theta x_0,T) \leq \phi(x_0,T)$,  we get from \eqref{eq:F-ep-xx-r} that, at $(x_0,T)$,
\[
\phi^2 +\left(\frac{1}{2}+\frac{1}{e^{x_0}-1} - \frac{x_0^2 e^{x_0}(e^{x_0}+1)}{2(e^{x_0}-1)^3} \right) \phi \leq 0
\]
This gives a contradiction as $\phi(x_0,T)^2>0$ and
\[
\frac{1}{2}+\frac{1}{e^{x_0}-1} - \frac{x_0^2 e^{x_0}(e^{x_0}+1)}{2(e^{x_0}-1)^3} \geq 0.
\]
Thus, we get 
\[
F^\ep_{xx} \leq 0 \qquad \text{ on } [0,\infty)^2.
\]
\end{proof}

By combining Lemmas \ref{lem:F-ep-increase}--\ref{lem:F-ep-concave}, we immediately get the following result.
\begin{lem}\label{lem:F-increase-and-concave}
    Assume $m\in (0,1]$ and \eqref{eq:rho-0}.
    Let $F$ be the unique solution to \eqref{eq:HJ}.
    Then, 
    \[
    0\leq F_x \leq m \qquad \text{ and } \qquad F_{xx} \leq 0 \qquad \text{ on } [0,\infty)^2
    \]
    in the viscosity sense.
\end{lem}

\subsection{Regularity of $F$ when $m\in\left(0,\frac{1}{2}\right)$}
In this section, we first prove the $C^{\infty}$ regularity of the solution $F$ to \eqref{eq:HJ} when $0<m<\frac{1}{2}$  in Proposition \ref{prop:Cinfty}.
Later, we obtain further regularity properties of $G$, which solves \eqref{eq:G}, in Proposition \ref{prop:derivatives-G}.

\smallskip

We start with the $C^{1,1}$ regularity result as follows.

\begin{prop}\label{prop:C1,1}
Assume the initial data $\rho_0\ :\ \mathbb{N}\to[0,\infty)$ satisfies $m_1(0)=m\in\left(0,\frac{1}{2}\right)$ and \eqref{eq:rho-0}. Let $F$ be the viscosity solution to \eqref{eq:HJ} with $F_0(x)=\sum_{j=1}^{\infty}(1-e^{-jx})\rho_0(j)$ for $x\geq0$. Then, $F\in C^{1,1}\left((0,\infty)^2\right)$.
\end{prop}

To prove this proposition, we introduce another viscosity term $a=a(x)$ given by a smooth, nondecreasing, and concave function such that
\begin{align}\label{def:a}
a(x)=\left\{\begin{array}{ll}
x,\hspace{0.5cm} x\in[0,1],\\
2,\hspace{5.3mm} x\in[3,\infty).
\end{array}\right.
\end{align}
For each $\varepsilon>0,$ let $F_2^{\varepsilon}$ be the viscosity solution to the equation \eqref{eq:add-ep} with $a$ defined as in \eqref{def:a}. We have that $F_2^{\varepsilon}\in C^{\infty}((0,\infty)^2)$. We first state relevant lemmas regarding the regularity of $F_2^{\varepsilon}$ on the boundary. 

\begin{lem}\label{lem:F_2ep}
Suppose $m\in\left(0,\frac{1}{2}\right)$ and  \eqref{eq:rho-0}. For each $\varepsilon>0$, let $F_2^{\varepsilon}$ be the solution to \eqref{eq:add-ep} with $a=a(x)$ defined as in \eqref{def:a}. Then, $F_2^{\varepsilon}$ is concave in $x$ and
\[
0\leq \partial_xF_2^{\varepsilon}\leq m\quad\text{ in }(0,\infty)^2.
\]
Moreover, $F_2^{\varepsilon}\in C^1([0,\infty)^2)$, and $\partial_xF_2^{\varepsilon}(0,t)=m$. Finally,
\[
\lim_{x\to0^+}x\partial^2_xF_2^{\varepsilon}(x,t)=0.
\]
\end{lem}

We skip the proof of this lemma, which is similar to those of \cite[Lemma 3.5, Lemma 3.6]{Tran-Van-1}.

\begin{lem}\label{lem:lowerboundofsecondderivative}
Suppose $m\in\left(0,\frac{1}{2}\right)$ and \eqref{eq:rho-0}. For each $\varepsilon>0$, let $F_2^{\varepsilon}$ be the solution to \eqref{eq:add-ep} with $a=a(x)$ defined as in \eqref{def:a}. Then, for $\varepsilon>0$ sufficiently small,
\begin{align*}
(1-e^{-x})\partial_x^2F_2^{\varepsilon}\geqslant -1\quad\text{ in }(0,\infty)^2.
\end{align*}
\end{lem}
\begin{proof}
The proof is inspired by that of \cite[Lemma 3.7]{Tran-Van-1} with various modifications. 
We break the proof into a few steps as follows.

\emph{Step 1}. Differentiating \eqref{eq:add-ep} twice in $x$, we get
\begin{multline}
\left(\partial_t\partial_x^2F_2^{\varepsilon}+\left[\partial_xF_2^{\varepsilon}-\left(m+\frac{1}{2}\right)\right]\partial_x^3F_2^{\varepsilon}\right)+(\partial_x^2F_2^{\varepsilon})^2\hfill\\
+\frac{1}{2}\partial_x^2F_2^{\varepsilon}+\frac{\partial_x^2F_2^{\varepsilon}}{e^x-1}-\frac{2e^x}{(e^x-1)^2}\partial_xF_2^{\varepsilon}+\frac{e^x(e^x+1)}{(e^x-1)^3}F_2^{\varepsilon}\\
\hfill=\varepsilon\left(a''\partial_x^2F_2^{\varepsilon}+2a'\partial_x^3F_2^{\varepsilon}+a\partial_x^4F_2^{\varepsilon}\right).\label{eq:diff-twice}
\end{multline}
Let
\[
Q^{\varepsilon}:=(1-e^{-x})\partial_x^2F_2^{\varepsilon}.
\]
Suppose that there exist $x_0,T>0$ such that
\[
\min_{[0,\infty)\times [0,T]}Q^{\varepsilon}(x,t)=Q^{\varepsilon}(x_0,T)=:\alpha<0.
\]
By the maximum principle, calculation shows that, at $(x_0,T)$,
\begin{align}
\begin{cases}
0\geq\partial_tQ^{\varepsilon}=(1-e^{-x_0})\partial_t\partial_x^2F_2^{\varepsilon},\\
0=\partial_xQ^{\varepsilon}\Longrightarrow\partial_x^2F_2^{\varepsilon}=(1-e^{x_0})\partial_x^3F_2^{\varepsilon},\label{eq:max-principle}\\
0\leq \partial_x^2Q^{\varepsilon}\Longrightarrow \partial_x^4F_2^{\varepsilon}\geq\frac{e^{x_0}+1}{(e^{x_0}-1)^2}\partial_x^2F_2^{\varepsilon}.
\end{cases}
\end{align}

\medskip

On one hand, we multiply the right-hand side of \eqref{eq:diff-twice} by $(1-e^{-x})^2$ and evaluate at $(x_0,T)$ to obtain a lower bound by
\begin{align}
&\varepsilon\left(a''(x_0)(1-e^{-x_0})^2\partial_x^2F_2^{\varepsilon}+2a'(x_0)(1-e^{-x_0})^2\partial_x^3F_2^{\varepsilon}+a(x_0)(1-e^{-x_0})^2\partial_x^4F_2^{\varepsilon}\right)\notag\\
&\geq\varepsilon\alpha\left(a''(x_0)(1-e^{-x_0})-2a'(x_0)e^{-x_0}+a(x_0)\frac{1+e^{-x_0}}{e^{x_0}-1}\right)\notag\\
&\geq Ce^{-x_0}\varepsilon\alpha\label{eq:lowerbound}
\end{align}
for some constant $C>0$ that depends only on $a=a(x)$. The last line follows from the choice of $a=a(x)$ such that $a''\leq0,a'\geq0$ and the function $\frac{a(x)(e^x+1)}{e^x-1}$ is bounded in $(0,\infty)$.
Also, \eqref{eq:max-principle} is repeatedly used here.

On the other hand, we multiply the left-hand side of \eqref{eq:diff-twice} by $(1-e^{-x})^2$ and evaluate at $(x_0,T),$ to obtain an upper bound by
\begin{align}
&\left((1-e^{-x_0})^2\partial_t\partial_x^2F_2^{\varepsilon}+\left[\partial_xF_2^{\varepsilon}-\left(m+\frac{1}{2}\right)\right](1-e^{-x_0})^2\partial_x^3F_2^{\varepsilon}\right)+((1-e^{-x_0})\partial_x^2F_2^{\varepsilon})^2\notag\\
+&\frac{e^{x_0}+1}{2(e^{x_0}-1)}(1-e^{-x_0})^2\partial_x^2F_2^{\varepsilon}-\frac{2e^{x_0}}{(e^{x_0}-1)^2}(1-e^{-x_0})^2\partial_xF_2^{\varepsilon}+\frac{e^{x_0}(e^{x_0}+1)}{(e^{x_0}-1)^3}(1-e^{-x_0})^2F_2^{\varepsilon}\notag\\
&\leq \alpha^2+e^{-x_0}\left(m+1-\partial_xF_2^{\varepsilon}+\frac{1}{2}e^{x_0}\right)\alpha+\frac{e^{-x_0}(e^{x_0}+1)}{e^{x_0}-1}F_2^{\varepsilon}-2e^{-x_0}\partial_xF_2^{\varepsilon}.\label{eq:upperbound}
\end{align}
Here, \eqref{eq:max-principle} is repeatedly used.

Connecting \eqref{eq:lowerbound} and \eqref{eq:upperbound}, we obtain
\[
\alpha^2+A\alpha+B\geq0,
\]
where
\begin{align*}
\begin{cases}
A:=e^{-x_0}\left(m+1-\partial_xF_2^{\varepsilon}\right)+\frac{1}{2}-Ce^{-x_0}\varepsilon\\
B:=\frac{e^{-x_0}(e^{x_0}+1)}{e^{x_0}-1}F_2^{\varepsilon}-2e^{-x_0}\partial_xF_2^{\varepsilon}.
\end{cases}
\end{align*}

\medskip

\emph{Step 2.} We first note that $F_2^{\varepsilon}(x,t)\leq m(1-e^{-x})$ by the maximum principle applied to \eqref{eq:add-ep}. We give an upper bound of $B$ by using this as\allowdisplaybreaks
\begin{align*}
B&=\frac{e^{-x_0}(e^{x_0}+1)}{e^{x_0}-1}F_2^{\varepsilon}-2e^{-x_0}\partial_xF_2^{\varepsilon}\\
&\leq\frac{e^{-x_0}(e^{x_0}+1)}{e^{x_0}-1}m(1-e^{-x_0})-2e^{-x_0}\partial_xF_2^{\varepsilon}\\
&=e^{-x_0}\left(2\delta+m(e^{-x_0}-1)\right)\\
&\leq \delta e^{-x_0}(1+e^{-x_0}),
\end{align*}
where $\delta:=m-\partial_xF_2^{\varepsilon}(x_0,T)$ satisfies $0\leq\delta\leq m$ by Lemma \ref{lem:F_2ep}. Consequently, for $0<m<\frac{1}{2}$ and $\varepsilon>0$ sufficiently small, we have 
\begin{align*}
A^2-4B&\geq\left(e^{-x_0}(\delta+1)+\frac{1}{2}-Ce^{-x_0}\varepsilon\right)^2-4\delta e^{-x_0}(1+e^{-x_0})\\
&\geq e^{-2x_0}(1-\delta)^2+\frac{1}{4}-e^{-x_0}(3\delta-1)-Ce^{-x_0}\varepsilon\\
&\geq \left(e^{-x_0}(1-\delta) -\frac{1}{2} \right)^2 + 2e^{-x_0}(1-2\delta- C\ep)>0.
\end{align*}

\medskip
By Lemma \ref{lem:low-up}, there exist two fixed constants $\ol h_- < \ol h_+ <-1/(2e)$ independent of $\ep$ and $x_0$ such that, for $\ep>0$ sufficiently small,
\[
\frac{-A-\sqrt{A^2-4B}}{2} < \ol h_- < \ol h_+ <\frac{-A+\sqrt{A^2-4B}}{2}.
\]
We now prove that
\[
\alpha\geq\frac{-A+\sqrt{A^2-4B}}{2}
\]
by verifying by contradiction that the other case
\[
\alpha\leq\frac{-A-\sqrt{A^2-4B}}{2}
\]
is impossible.

Assume for the contrary that the latter holds for some fixed small $\varepsilon>0$. From the form $F_0(x)=\sum_{j=1}^{\infty}\left(1-e^{-jx}\right)\rho_0(j)$, we have that
\[
\ol h_+ <-\frac{1}{2e}\leq-\frac{m}{e}\leqslant (1-e^{-x})F_0''(x)\leqslant0.
\]
By continuity, there exists $T^{\varepsilon}\in(0,T)$ so that
\[
\ol h_-<\min_{(0,\infty)\times[0,T^{\varepsilon}]}Q^{\varepsilon}<\ol h_+,
\]
which is a contradiction.
\end{proof}

\begin{lem}\label{lem:low-up}
Fix $\delta \in [0,1/2]$.
For $r\in [0,1]$, denote by
\begin{align*}
    h_-(r)&=-\frac{1}{2} - r(1+\delta) -\sqrt{\frac{1}{4}+r^2(1-\delta)^2 + r(1-3\delta)}\ ,\\
    h_+(r)&=-\frac{1}{2} - r(1+\delta) +\sqrt{\frac{1}{4}+r^2(1-\delta)^2 + r(1-3\delta)}\ .
\end{align*}
Then,
\[
\max_{[0,1]} h_- = -1 < \min_{[0,1]} h_+.
\]
\end{lem}

\begin{proof}
It is straightforward that $h_-$ is a nonincreasing function and
\[
\max_{[0,1]} h_- = h_-(0)=-1.
\]
By a direct comparison, we have that, for $r\in [0,1]$,
\[
h_+(r) > -1,
\]
which concludes the proof.
\end{proof}

\begin{proof}[Proof of Proposition \ref{prop:C1,1}]
The $C^{1,1}$ regularity of $F$ follows immediately from the results in Lemmas \ref{lem:F_2ep} and \ref{lem:lowerboundofsecondderivative}.
\end{proof}

\medskip

Here is the main result of this subsection.

\begin{prop}\label{prop:Cinfty}
    Under the same setting of Proposition \ref{prop:C1,1} with the same notations, it holds that $F\in C^{\infty}\left((0,\infty)^2\right)$. Moreover, for each $k\in\mathbb{N}$ and $R>0$, there exists a constant $C=C(k,R)>0$ such that
    \begin{equation}\label{eq:bound-F-k}
    \left\|x^{k-1} \partial^k_x F(x,t)\right\|_{L^\infty((0,R)^2)} \leq C(k,R).
    \end{equation}
\end{prop}

We omit the proof of Proposition \ref{prop:Cinfty} as it can be verified similarly to \cite[Proposition 3.9, Lemma 3.11, Lemma 3.12]{Tran-Van-1}.
Basically, once we have $F\in C^{1,1}((0,\infty)^2)$, we use the method of characteristics to get $F\in C^{\infty}\left((0,\infty)^2\right)$.
Then, we differentiate \eqref{eq:CF} with respect to $x$ for $k$ times and use the localization around the characteristics and the maximum principle to get \eqref{eq:bound-F-k}.

We emphasize that we are not able to prove directly that $F$ is a Bernstein function, that is, $(-1)^{k+1}\partial^k_x F \geq 0$ in $(0,\infty)^2$ for all $k\in \N$.
Technically, the zeroth order term $\frac{1}{2}F + \frac{1}{e^x-1}F$ in \eqref{eq:HJ} becomes extremely convoluted after repeated differentiations in $x$ and we are not able to handle it.
Instead, we resort our approach to considering \eqref{eq:G} in the next section.
The results in Section \ref{sec:solution-CF} are new and not analogous to any of the results in \cite{Tran-Van-1}.

\section{Regularity of solutions to \eqref{eq:G} and existence of mass-conserving solutions to \eqref{eq:CF} when $m\in (0,1/2)$}\label{sec:solution-CF}

Here, we work with $G$ instead of $F$. 
Recall that
\[
G(z,t)=F(-\log z,t) \qquad \text{ for } (z,t)\in (0,1]\times [0,\infty).
\]
Then, \eqref{eq:HJ} can be formulated in an equivalent way as
\begin{equation}\label{eq:G-rep}
    \begin{cases}
        G_t +\frac{1}{2}\left(zG_z+m\right)\left(zG_z+m+1\right)+ \frac{1+z}{2(1-z)}G -m=0 \qquad &\text { in } (0,1)\times(0,\infty),\\
        0\leq G(z,t) \leq m(1-z) \qquad &\text { in } (0,1)\times(0,\infty),\\
        G(z,0)=G_0(z) \qquad &\text { on } (0,1].
    \end{cases}
\end{equation}
Here, $G_0(z)=F_0(-\log z)$ for $z\in (0,1]$.

Thanks to Proposition \ref{prop:Cinfty} and the elementary fact that $0\leq 1-z\leq-\log z$ for $z\in(0,1]$, for each $k\in\mathbb{N}$ and $R>0$, we have
\begin{align}\label{eq:boundary-G}
\left\|(1- z)^{k-1} \partial^k_z G(z,t)\right\|_{L^\infty((1/2,1)\times (0,R))}\leq C(k,R).
\end{align}
for some constant $C(k,R)>0$.

\begin{prop}\label{prop:derivatives-G}
    For $k\in \N$,
    \begin{equation}\label{eq:derivative-k-G}
    \partial^k_z G \leq 0 \qquad \text{ in } (0,1)\times(0,\infty).
    \end{equation}
\end{prop}
\begin{proof}
    We prove \eqref{eq:derivative-k-G} by induction.
    Thanks to Lemma \ref{lem:F-ep-increase}, \eqref{eq:derivative-k-G} holds when $k=1$.

    \smallskip

    Assume that \eqref{eq:derivative-k-G} is true for $k\leq n-1$ for a given natural number $n\geq 2$.
    We now show that \eqref{eq:derivative-k-G} also holds for $k=n$.
    We give a proof of this point by contradiction.

    \smallskip

    First, we can rewrite \eqref{eq:G-rep} as
    \[
    G_t + \frac{1}{2}z^2 (G_z)^2 + \left(m+\frac{1}{2}\right) zG_z -\frac{1}{2}G +\frac{G}{1-z}+\frac{1}{2}m(m-1)=0.
    \]
    We differentiate this equation $n$ times with respect to $z$ and aim at writing an equation for $\varphi(z,t)=(1-z)^{n-1} \partial^n_z G(z,t)$.
    To do this, we compute
    \begin{align*}
        \frac{G(z,t)}{1-z} &= \frac{G(z,t)-G(1,t)}{1-z}\\
        &=-\int_0^1 G_z(1+r(z-1),t)\,dr.
    \end{align*}
    Hence,
    \begin{align*}
        &(1-z)^{n-1}\partial_z^n\left(\frac{G(z,t)}{1-z}\right) \\
        =\ & -(1-z)^{n-1}\int_0^1 \partial_z^{n+1}G(1+r(z-1),t)r^n\,dr\\
        =\ & \frac{\varphi(z,t)}{1-z} - \frac{n}{1-z}\int_0^1 \varphi(1+r(z-1),t)\,dr.
    \end{align*}
    We compute next
    \begin{align*}
        (1-z)^{n-1}\partial^n_z(z G_z)&=(1-z)^{n-1}\left(n \partial^n_z G + z \partial^{n+1}_z G\right)\\
        &=n \varphi + z\left(\varphi_z + \frac{n-1}{1-z}\varphi\right).
    \end{align*}
    Finally,
    \begin{align*}
        &(1-z)^{n-1}\partial^n_z(z^2 (G_z)^2)\\
       =\ & (1-z)^{n-1} \left(z^2 \partial^n_z((G_z)^2) + 2nz \partial^{n-1}_z((G_z)^2) + n(n-1) \partial^{n-2}_z ((G_z)^2)\right)\\
       =\ &2z^2 G_z \left(\varphi_z + \frac{n-1}{1-z}\varphi\right)+2nz^2 G_{zz} \varphi +4nzG_z \varphi  + h(z,t),
    \end{align*}
    where
    \begin{multline*}
        h(z,t)=(1-z)^{n-1} \Bigg[ z^2 \sum_{i=2}^{n-2} \binom{n}{i} \partial^{i+1}_z G\cdot  \partial^{n+1-i}_z G \\
        +2nz \sum_{i=1}^{n-2}\binom{n-1}{i} \partial^{i+1}_z G\cdot  \partial^{n-i}_z G + n(n-1)\sum_{i=0}^{n-2}\binom{n-2}{i} \partial^{i+1}_z G\cdot  \partial^{n-1-i}_z G
        \Bigg] \geq 0
    \end{multline*}
    in light of the induction hypothesis.

    \smallskip

    Combine everything to yield
    \begin{equation}\label{eq:varphi}
    \varphi_t + \left(zG_z+m+\frac{1}{2}\right) z\varphi_z + q(z,t) \varphi - \frac{n}{1-z}\int_0^1 \varphi(1+r(z-1),t)\,dr + \frac{1}{2}h(z,t)=0,    
    \end{equation}
    where
    \[
    q(z,t)=nz^2G_{zz}+\left(\frac{2n-(n+1)z}{1-z}\right)zG_z+\left(m+\frac{1}{2}\right)\left(\frac{n-z}{1-z}\right)+\frac{1+z}{2(1-z)}.
    \]

    By the localization argument, we may assume without loss of generality that there exist $\delta\in (0,1)$ and $T>0$ such that
    \[
    \begin{cases}
    \xi(t)=\varphi(z_t,t) = \max_{z\in(0,1]} \varphi(z,t) \quad &\text{ for } t\in [0,T],\\
    z_t \leq 1-\delta \quad &\text{ for } t \in [0,T],\\
    \max_{t\in [0,T]} \xi(t)>0.
    \end{cases}
    \]
    Use this and the maximum principle for \eqref{eq:varphi} to deduce
    \[
    \xi'(t) +\left( q(z_t,t) - \frac{n}{1-z_t}\right) \xi(t) \leq 0.
    \]
    By the usual differential inequality, we imply that $\xi \leq 0$, which is absurd.
    The proof is complete.
\end{proof}

\medskip

The following shows we can extract the densities of size $l\in\mathbb{N}$.

\begin{lem}\label{lem:rho-l}
For each $l\in\mathbb{N}$, there exists a $\rho(l,t)\geq0$ such that
\begin{align}\label{eq:G-series}
G(z,t)=\sum_{l=1}^{\infty}\left(1-z^l\right)\rho(l,t)\quad\text{ for }(z,t)\in[0,1]\times(0,\infty).
\end{align}
Moreover, we have
\begin{align}\label{eq:m0-m1}
G(0,t)=\sum_{l=1}^{\infty}\rho(l,t)\quad\text{ and }\quad m=\sum_{l=1}^{\infty}l\rho(l,t).
\end{align}
\end{lem}
\begin{proof}
For each $l\in\mathbb{N},\ t>0$, let
\begin{align}
\rho(l,t):=-\frac{1}{l!}\partial_z^lG(0,t),\label{eq:choice-rho-l}
\end{align}
which is understood as the right-continuous extension of $-\frac{1}{l!}\partial_z^lG(z,t)$. Note that this is well-defined and nonnegative thanks to Lemma \ref{prop:derivatives-G}.

By Taylor's expansion in integral form, for each $0\leq a<z<1,\ t>0,\ k\in\mathbb{N}$, we have
\begin{align}
G(z,t)=G(a,t)+\sum_{l=1}^{k}\frac{\partial_z^lG(a,t)}{l!}(z-a)^l+R_k(z,t;a)\label{eq:G-taylor-expn}
\end{align}
where
\begin{align}
R_k(z,t;a):=\int_a^z\frac{\partial_z^{k+1}G(x,t)}{k!}(z-x)^kdx.\label{eq:remainder-term}
\end{align}
As $R_k\leq0$, we have
\[
0\leq-\sum_{l=1}^{k}\frac{\partial_z^lG(a,t)}{l!}(z-a)^l\leq G(a,t)-G(z,t),
\]
and thus, we have the convergence of the series $-\sum_{l=1}^{\infty}\frac{\partial_z^lG(a,t)}{l!}(z-a)^l$. In particular, it holds that
\[
-\frac{\partial_z^kG(a,t)}{k!}(z-a)^k\leq G(a,t)-G(z,t)\leq 1.
\]

\medskip

We note that for each $0\leq a<x<z<1,\ t>0,$
\[
-\frac{\partial_z^{k+1}G(x,t)}{k!}(z-x)^k=-\frac{\partial_z^{k+1}G(x,t)}{(k+1)!}(z_0-x)^{k+1}\cdot\frac{k+1}{z_0-x}\cdot\left(\frac{z-x}{z_0-x}\right)^{k}\leq\frac{k+1}{z_0-z}\cdot\left(\frac{z}{z_0}\right)^{k},
\]
where $z_0$ can be taken any number in $(z,1)$. Therefore, we have for each $0\leq a<x<z<1,\ t>0$,
\begin{align*}
0\leq-R_k(z,t;a)\leq\frac{k+1}{z_0-z}\cdot\left(\frac{z}{z_0}\right)^{k},
\end{align*}
where the right-hand side converges to 0 as $k\to\infty$. Therefore,
\begin{equation}\label{eq:G-series-new}
G(z,t)=G(0,t)-\sum_{l=1}^{\infty}\rho(l,t)z^l\quad\text{ for }z\in[0,1].
\end{equation}
The endpoint $z=1$ can be included by applying Abel's Lemma. Moreover, evaluating at $z=1$ gives 
\[
G(0,t)=\sum_{l=1}^{\infty}\rho(l,t).
\]
Differentiate \eqref{eq:G-series-new} in $z$, apply Abel's Lemma once more, and use the fact that $G_z(1,t)=m$ to obtain 
\[
m=\sum_{l=1}^{\infty}l\rho(l,t).
\]
The proof is complete.
\end{proof}

We next state and prove the regularity of $G$ on the boundary $z=0$. Let $\rho(0,t)=-G(0,t)$ for convenience.

\begin{lem}\label{lem:G-boundary-regularity}
For each $l\in\mathbb{N}\cup\{0\}$, the function $\rho(l,t)$ is differentiable in $t>0$. 
\end{lem}
\begin{proof}
Let $l\in\mathbb{N}\cup\{0\}$. Differentiate the first equation of \eqref{eq:G-rep} $l$ times in $z$ to obtain
\[
\partial_t\partial_z^lG=-\frac{1}{2}\partial_z^l\left((zG_z)^2\right)-\left(m+\frac{1}{2}\right)\partial_z^l\left(zG_z\right)-\partial_z^l\left(\frac{1+z}{2(1-z)}G\right).
\]
Fix $t>0$. For each $s\in (0,\infty)\setminus \{t\}$, by the mean value theorem, there exists $r=r(s)>0$ between $s$ and $t$  such that
\begin{align}
&\partial_t\rho(l,t)=-\frac{1}{l!}\partial_t\left(\partial_z^lG(0,t)\right)\notag\\
=\ &-\frac{1}{l!}\lim_{s\to t}\frac{\partial_z^lG(0,s)-\partial_z^lG(0,t)}{s-t}\notag\\
=\ &-\frac{1}{l!}\lim_{s\to t}\lim_{z\to0^+}\frac{1}{s-t}\left(\partial_z^lG(z,s)-\partial_z^lG(z,t)\right)\notag\\
=\ &-\frac{1}{l!}\lim_{s\to t}\lim_{z\to0^+}\partial_t\partial_z^lG(z,r)\notag\\
=\ &-\frac{1}{l!}\lim_{s\to t}\lim_{z\to0^+}\left(-\frac{1}{2}\partial_z^l\left((zG_z)^2\right)-\left(m+\frac{1}{2}\right)\partial_z^l\left(zG_z\right)-\partial_z^l\left(\frac{1+z}{2(1-z)}G\right)\right)\biggr\rvert_{(z,r)}\notag\\
=\ &-\frac{1}{l!}\left(-\frac{1}{2}\partial_z^l\left((zG_z)^2\right)-\left(m+\frac{1}{2}\right)\partial_z^l\left(zG_z\right)-\partial_z^l\left(\frac{1+z}{2(1-z)}G\right)\right)\biggr\rvert_{(0,t)}.\label{eq:rho_t-of-l}
\end{align}
In the last line, we used the continuity of $G,\cdots,\partial_z^{l+1}G$ in $(z,t)\in[0,1)\times(0,\infty)$. This can be seen from the fact that $\|\partial_z^kG\|_{L^{\infty}\left(\left[0,\frac{1}{2}\right)\times(0,\infty)\right)}\leq C(k)$ for each $k\in\mathbb{N}\cup\{0\}$, which follows from \eqref{eq:G-series}.
\end{proof}

\medskip

\begin{proof}[Proof of Theorem \ref{thm:existence-of-C-F}]
Using \eqref{eq:G-series} and \eqref{eq:m0-m1}, we can check that our choice \eqref{eq:choice-rho-l} of the densities solves the equation \eqref{eq:CF} with the kernels as in \eqref{eq:ab} by a simple substitution of \eqref{eq:G-series} into \eqref{eq:rho_t-of-l}.
\end{proof}

\section{Equilibria and large time behaviors}\label{sec:large-time-behavior}

\subsection{Large time behaviors}
We prove Theorem \ref{thm:m3} in this section.
We start with the convergence of $G(z,t)$ as $t\to\infty$.

\begin{lem}\label{lem:G-large-time}
Assume the initial data $\rho_0 : \mathbb{N}\to[0,\infty)$ has mass $m\in (0,1]$ and satisfies \eqref{eq:rho-0}.
The solution $G(\cdot,t)$ to \eqref{eq:G} converges locally uniformly on $(0,1]$ to $\tilde G$ as $t\to \infty$, which is the unique viscosity solution  to
    \begin{equation}\label{eq:G-stationary}
    \begin{cases}
        \frac{1}{2}\left(z\tilde G_z+m\right)\left(z\tilde G_z+m+1\right)+ \frac{1+z}{2(1-z)}\tilde G -m=0 \qquad &\text { in } (0,1),\\
        0\leq\tilde G(z)\leq m(1-z) \qquad &\text { on } [0,1],
    \end{cases}
    \end{equation}
    as $t\to\infty$.
\end{lem}

\begin{proof}
We prove the corresponding result for the solution $F(x,t)$ to \eqref{eq:CF} for convenience reason as the cut-off function $\theta_n(x)=\max\left\{\frac{1}{n},e^x-1\right\}$ for $x\in[0,\infty)$ is already introduced in Section 2.

From the proofs of Lemma \ref{lem:F-n} and Proposition \ref{prop:m1}, we have \[
0\leq F(x,t) \leq m(1-e^{-x})
\]
for $(x,t)\in[0,\infty)^2.$ Accordingly, the envelopes 
\[
\overline{F}(x):=\limsup_{t\to\infty}\!{}^{*}F(x,t),\qquad \underline{F}(x):=\liminf_{t\to\infty}\!{}_*\, F(x,t)
\]
satisfy
\[
0\leq \underline{F}(x) \leq \overline{F}(x) \leq m(1-e^{-x})
\]
for $x\in[0,\infty)$. We note that the function $\underline{F}(x)+\frac{1}{n}$ is a viscosity supersolution to
\begin{equation*}
    \begin{cases}
        \frac{1}{2}\left(\tilde F_x-m\right)\left(\tilde F_x-m-1\right)+\frac{1}{2}\tilde F + \frac{\tilde F}{  \theta_n} -m=0 \qquad &\text { in } (0,\infty),\\
        \tilde F(0)=0,
    \end{cases}
    \end{equation*}
while the function $\overline{F}(x)$ is a viscosity subsolution to the above. Thanks to the monotone term $\frac{1}{2}F$, the comparison principle for bounded solutions is applicable, which allows us to deduce that $\overline{F}(x)\leq\underline{F}(x)+\frac{1}{n}$ for every $n\in\N$. Letting $n\to\infty$ completes the proof.
\end{proof}

\begin{cor}\label{cor:m=1}
    Assume the initial data $\rho_0 : \mathbb{N}\to[0,\infty)$ has mass $m=1$ and satisfies \eqref{eq:rho-0}.
The solution $G(\cdot,t)$ to \eqref{eq:G} converges locally uniformly on $(0,1]$ to $\tilde G \equiv 0$ as $t\to \infty$.
\end{cor}

\begin{proof}
    Thanks to Lemma \ref{lem:G-large-time}, $G(\cdot,t)$ to \eqref{eq:G} converges locally uniformly on $(0,1]$ to $\tilde G$, which is the viscosity solution to \eqref{eq:G-stationary}.
    For $m=1$, it is clear that $\tilde G \equiv 0$ is the unique viscosity solution to \eqref{eq:G-stationary}.
\end{proof}

\medskip

\begin{rem}\label{rem:m=1}
    We consider the initial data $\rho_0 : \mathbb{N}\to[0,\infty)$ has mass $m=1$ and satisfies \eqref{eq:rho-0}.
    By Corollary \ref{cor:m=1}, both $F(\cdot,t)$ and $G(\cdot,t)$ converge locally uniformly to $0$ as $t\to \infty$.
    This implies that gelation happens and all the mass is gone at infinity, which agrees with the heuristic discussion in Section \ref{sec:formal cal}.
    The situation here is completely different from the continuous case where there is an one-parameter family of stationary solutions \cite{Mitake-Tran-Van}.
\end{rem}

\medskip

\begin{lem}\label{lem:rho-l-large-time}
Assume the initial data $\rho_0 : \mathbb{N}\to[0,\infty)$ has mass $m\in\left(0,\frac{1}{2}\right)$ and satisfies \eqref{eq:rho-0}.
The density $\rho(l,t)$ defined as in \eqref{eq:choice-rho-l} converges as $t\to\infty$ for each $l\in\mathbb{N}$.
\end{lem}
\begin{proof}
For $l\in\mathbb{N}$ and $h\in\left(0,\frac{1}{2l}\right)$, write the $l$-th order difference of $G$ at $z=0$ as
\[
\Delta_h^lG(0,t):=\sum_{j=0}^l(-1)^{l-j}\binom{l}{j}G(jh,t).
\]
By a successive use of the mean value theorem, there exists $\theta\in(0,lh)$ such that
\[
\Delta_h^lG(0,t)=\partial_z^lG(\theta,t).
\]
Therefore,
\begin{align*}
\left|\Delta_h^lG(0,t)-\partial_z^lG(0,t)\right|=\left|\partial_z^lG(\theta,t)-\partial_z^lG(0,t)\right|\leq C_l h,
\end{align*}
where $C_l:=l\cdot\sup_{0\leq z\leq\frac{1}{2}}\left|\partial_z^{l+1}G(z,t)\right|\leq 2^{l+2}(l+1)!l$.
Letting $t\to\infty$ in the above to yield
\begin{align*}
-C_lh\leq\Delta_h^l\tilde G(0)-\limsup_{t\to\infty}\partial_z^lG(0,t)\leq\Delta_h^l\tilde G(0)-\liminf_{t\to\infty}\partial_z^lG(0,t)\leq C_lh.
\end{align*}
Here, $\Delta_h^l\tilde G(0)=\sum_{j=0}^l(-1)^{l-j}\binom{l}{j}\tilde G(jh)$.
We hence get
\[
0\leq \limsup_{t\to\infty}\partial_z^lG(0,t)-\liminf_{t\to\infty}\partial_z^lG(0,t)\leq 2C_lh.
\]
 We finish the proof by letting $h\to0$.
\end{proof}

\medskip

\begin{lem}\label{lem:tilde-G-series}
Assume the initial data $\rho_0 : \mathbb{N}\to[0,\infty)$ has mass $m\in\left(0,\frac{1}{2}\right)$ and satisfies \eqref{eq:rho-0}.
Let $\tilde{\rho}(l)=\lim_{t\to\infty}\rho(l,t)$ for each $l\in\mathbb{N}$. Then, the solution $\tilde G$ to \eqref{eq:G-stationary} satisfies
\begin{align}
\tilde G(z)=\sum_{l=1}^{\infty}\left(1-z^{l}\right)\tilde{\rho}(l)\quad\text{ for }z\in[0,1].\label{eq:G-stationary-series}
\end{align}
In particular, we have
\begin{align}\label{eq:m0-m1-stationary}
m(1-m)=\sum_{l=1}^{\infty}\tilde{\rho}(l)\quad\text{ and }\quad m=\sum_{l=1}^{\infty}l\tilde{\rho}(l).
\end{align}
\end{lem}
\begin{proof}
We start with an estimate of the remainder term $R_k(z,t;a)$ (with $a=0$) defined as in \eqref{eq:remainder-term} by
\begin{align*}
\left|R_k(z,t;0)\right|&\leq\int_0^z\frac{\left|\partial_z^{k+1}G(x,t)\right|}{k!}(z-x)^kdx\\
&\leq\int_0^z\frac{(k+1)!}{(1-x)^{k+2}}\cdot\frac{(z-x)^k}{k!}dx\leq\frac{(k+1)z^{k+1}}{(1-z)^2}
\end{align*}
for a fixed $z\in[0,1)$.

Taking $t\to\infty$ in \eqref{eq:G-taylor-expn} with $a=0$, we have
\[
\left|\tilde G(z)-\left(\tilde G(0)-\sum_{l=1}^k\tilde{\rho}(l)z^l\right)\right|\leq\frac{(k+1)z^{k+1}}{(1-z)^2}.
\]
for each $z\in[0,1)$. Letting $k\to\infty$ gives \eqref{eq:G-stationary-series} using the fact that $\tilde G(1)=0$. Also, \eqref{eq:m0-m1-stationary} follows immediately.
\end{proof}

\medskip

\begin{lem}\label{lem:recursive-formula}
Let $\tilde{\rho}(l)$ be as in Lemma \ref{lem:tilde-G-series} for each $l\in\mathbb{N}$. Then, the sequence $\{\tilde{\rho}(l)\}_{l\in\mathbb{N}}$ satisfies the recursive formula
\begin{align*}
\tilde{\rho}(l)=\frac{2m(1-m)}{(2m+1)l+1}+ \frac{1}{(2m+1)l+1} \left(\sum_{i=1}^{l-1}i(l-i) \tilde{\rho}(i)\tilde{\rho}(l-i)-2 \sum_{i=1}^{l-1}\tilde{\rho}(i) \right).   
\end{align*}
\end{lem}
\begin{proof}
We derive by substituting \eqref{eq:G-stationary-series} into \eqref{eq:G-stationary} that
\begin{equation}\label{eq:coefficient-vanish}
\sum_{l=1}^{\infty}\left(\sum_{i=l}^{\infty}\tilde{\rho}(i)+\frac{1}{2}\sum_{i=1}^{l-1}i(l-i)\tilde{\rho}(i)\tilde{\rho}(l-i)-\left(\left(m+\frac{1}{2}\right)l+\frac{1}{2}\right)\tilde{\rho}(l)\right)z^l=0,
\end{equation}
where the sum $\sum_{i=1}^{l-1}i(l-i)\tilde{\rho}(i)\tilde{\rho}(l-i)$ is understood as zero for $l=1$. As we have $\sum_{i=l}^{\infty}\tilde{\rho}(i)=m(1-m)-\sum_{i=1}^{l-1}\tilde{\rho}(i)$, and as the coefficient of $z^l$ in \eqref{eq:coefficient-vanish} must vanish for each $l\in\mathbb{N}$, we conclude that the recursive formula holds for each $l\in\mathbb{N}$.
\end{proof}

\medskip

\begin{proof}[Proof of Theorem \ref{thm:m3}]
The theorem follows immediately from Lemmas \ref{lem:rho-l-large-time}, \ref{lem:tilde-G-series}, \ref{lem:recursive-formula}.
\end{proof}

\subsection{Equilibria}

\begin{lem}\label{lem:series-m-big}
    Fix $m \geq 1$.
    Then, \eqref{eq:CF} with kernels given as in \eqref{eq:ab} does not admit a stationary solution $\tilde \rho:\N \to [0,\infty)$ such that
    \begin{equation*}
    \sum_{l=1}^\infty l \tilde \rho(l) =m.
    \end{equation*}
\end{lem}

\begin{proof}
    This is an immediate consequence of Lemma \ref{lem:recursive-formula}. 
    Assume otherwise that \eqref{eq:CF} with kernels given as in \eqref{eq:ab} admits a stationary solution $\tilde \rho:\N \to [0,\infty)$ satisfying
    \begin{equation*}
    \sum_{l=1}^\infty l \tilde \rho(l) =m.
    \end{equation*}
Thanks to Lemma \ref{lem:recursive-formula},
\[
\tilde{\rho}(1)=\frac{m(1-m)}{m+1},
\]
which is negative when $m>1$, which gives a contradiction. 
When $m=1$, we have $\tilde{\rho}(l)=0$ for all $l\in\mathbb{N}$, which also gives a contradiction.
\end{proof}

We note that in the above lemma, the case $m=1$ agrees with Corollary \ref{cor:m=1} and Remark \ref{rem:m=1}, and the case $m>1$ is consistent with the result of Theorem \ref{thm:nonexistence-of-C-F}.

\medskip

Thanks to Lemma \ref{lem:tilde-G-series} and \ref{lem:recursive-formula} and a continuity argument, we have the following result.
Here, the end point $m=1/2$ is included because of the continuity with respect to $m$ of $\tilde \rho(l)=\tilde \rho_m(l)$ in \eqref{eq:recursive}.

\begin{lem}\label{lem:series-m-0-1/2}
    Fix $m \in (0,1/2]$.
    Then, \eqref{eq:CF} with kernels given as in \eqref{eq:ab} admits a unique stationary solution $\tilde \rho:\N \to [0,\infty)$ such that $m_1(\tilde{\rho})=m$. Further, for $l\in \N$,
\begin{align*}
\tilde{\rho}(l)=\frac{2m(1-m)}{(2m+1)l+1}+ \frac{1}{(2m+1)l+1} \left(\sum_{i=1}^{l-1}i(l-i) \tilde{\rho}(i)\tilde{\rho}(l-i)-2 \sum_{i=1}^{l-1}\tilde{\rho}(i) \right).   
\end{align*}
\end{lem}



\end{document}